\documentclass{amsart}
\usepackage{amssymb}
\usepackage{amsmath}

\newtheorem{theorem}{Theorem}[section]
\newtheorem{lemma}[theorem]{Lemma}
\newtheorem{proposition}[theorem]{Proposition}
\theoremstyle{remark}
\newtheorem*{remark}{Remark}
\newtheorem*{remarks}{Remarks}
\numberwithin{equation}{section}
\begin{document}

\title[Relations for annihilating fields]{Relations for annihilating fields of
standard modules for affine Lie algebras}

\author{Mirko Primc}
\address{University of Zagreb\\
Department of Mathematics\\
Bijeni\v{c}ka 30, Zagreb, Croatia} \email{primc@math.hr}
\thanks{Partially supported by the Ministry of Science
and Technology of the Republic of Croatia, grant 037002.}
\begin{abstract}
J.~Lepowsky and R.~L.~Wilson initiated the approach to
combinatorial Rogers-Ramanujan type identities via the vertex
operator constructions of representations of affine Lie algebras.
In a joint work with Arne Meurman this approach is developed
further in the framework of vertex operator algebras. The main
ingredients of that construction are defining relations for
standard modules and  relations among them. The arguments involve
both representation theory and combinatorics, the final results
hold only for affine Lie algebras $A_1^{(1)}$ and $A_2^{(1)}$. In
the present paper some of those arguments are formulated and
extended for general affine Lie algebras. The main result is a
kind of rank theorem, guaranteeing the existence of combinatorial
relations among relations, provided that certain purely
combinatorial quantities are equal to dimensions of certain
representation spaces. Although the result holds in quite general
setting, applications are expected mainly for standard modules of
affine Lie algebras.
\end{abstract}
\maketitle

\section{Introduction}

J.~Lepowsky and R.~L.~Wilson gave in \cite{LW} a Lie-theoretic
interpretation and proof of the classical Rogers-Rama\-nu\-jan
identities in terms of representations of the affine Lie algebra
$\tilde{\mathfrak g}={\mathfrak sl}(2,\mathbb C)\,\widetilde{}$\,.
The identities are obtained by expressing in two ways the
principal characters of vacuum spaces  for the principal
Heisenberg subalgebra of $\tilde{\mathfrak g}$. The product sides
follow {}from the principally specialized Weyl-Kac character
formula for level 3 standard $\tilde{\mathfrak g}$-modules; the
sum sides follow {}from the vertex operator construction of bases
of level 3 standard $\tilde{\mathfrak g}$-modules, parametrized by
partitions satisfying difference 2 conditions. In fact, Lepowsky
and Wilson gave a construction of combinatorial bases of all
standard $\tilde{\mathfrak g}$-modules: the spanning follows in
general {}from the vertex operator ``generalized anticommutation
relations'', but the Lie-theoretic proof of linear independence
given for level 3 modules could not be extended for higher levels
$k\geq 4$. In \cite{MP1} a Lie-theoretic proof of linear
independence of these bases was given. To be precise, the vertex
operator construction of combinatorial bases of maximal submodules
of the corresponding Verma modules was given, so that the linear
independence was easy, and the spanning followed {}from the
relations among the generalized anticommutation relations. Similar
ideas have been used in \cite{FNO} and \cite{C}.

In \cite{MP2} these ideas were applied for the affine Lie algebra
$\tilde{\mathfrak g}={\mathfrak sl}(2,\mathbb C)\,\widetilde{}$\,
in the homogeneous picture by using the methods of vertex operator
algebras (cf. \cite{FLM}, \cite{FHL}). In this setting the
``generalized anticommutation relations'' $x_\theta(z)^{k+1}=0$,
studied before in \cite{LP}, are seen as annihilating fields of
standard modules. For level $k$ standard $\tilde{\mathfrak
g}$-modules the annihilating fields are associated with vectors in
the maximal ideal $N^1(k\Lambda_0)$ in the universal vertex
operator algebra $N(k\Lambda_0)$. The maximal ideal is generated
by a finite dimensional irreducible ${\mathfrak g}$-module $R$,
and the coefficients of annihilating fields associated to elements
in $R$ form a loop $\tilde{\mathfrak g}$-module $\bar R$. The main
point is that $\bar R$ is a set of defining relations for standard
$\tilde{\mathfrak g}$-modules: for every Verma module $M(\Lambda)$
with integral dominant $\Lambda$ the maximal submodule
$M^1(\Lambda)$ can be written as
$$
M^1(\Lambda)=\bar R M(\Lambda),\qquad L(\Lambda)=M(\Lambda)/\bar R M(\Lambda).
$$

 Since the character of Verma module $M(\Lambda)$ is easily
described, a combinatorial description of the character of
standard module $L(\Lambda)$ may be obtained {}from a
combinatorial description of $M^1(\Lambda)$. So in
Lepowsky-Wilson's approach a construction of combinatorial basis
of $L(\Lambda)$ can be obtained by constructing first a
combinatorial basis of $M^1(\Lambda)$. If $v_\Lambda\in
M(\Lambda)$ is a highest weight vector, then it is natural to seek
for a combinatorial basis of $M^1(\Lambda)=\bar R M(\Lambda)$
within the spanning set
\begin{equation}\label{1.1}
r\,x_1x_2\dots x_s v_\Lambda,\qquad r\in \bar R, \ x_i\in
\tilde{\mathfrak g}.
\end{equation}

It is clear that it is enough to consider elements $x_i\in
\tilde{\mathfrak g}$ {}from some fixed basis $\tilde B$ in
$\tilde{\mathfrak g}$. The corresponding monomial basis $\mathcal
P$ of the symmetric  algebra $S(\tilde{\mathfrak g})$ can be used
for parametrization of a monomial basis of the universal
enveloping algebra $U(\tilde{\mathfrak g})$: a monomial basis
element $u(\pi)=x_1x_2\dots x_s\in U(\tilde{\mathfrak g})$
corresponds to $\pi=x_1x_2\dots x_s\in\mathcal P$. The elements
$\pi\in\mathcal P$ are interpreted as colored partitions.

In the same manner, there is a basis of the space of relations
$\bar R$ parametrized by elements in $\mathcal P$. For a fixed
order $\preceq $ on $\mathcal P$ every nonzero element $r\in\bar
R$ is expanded in terms of monomial basis elements
$$
r=c_\rho u(\rho)+\sum_{\pi\succ\rho}c_\pi u(\pi),
$$
where $c_\rho, c_\pi\in\mathbb C$, $c_\rho\neq 0$. The first
nonzero monomial $u(\rho)$ which appears in this expansion is
called the leading term and denoted as $\rho=\ell \!\text{{\it
t\,}}(r)$. Then there is a basis $\{r(\rho)\mid\rho\in\ell
\!\text{{\it t\,}}(\bar R)\}$ of $\bar R$ such that $\ell
\!\text{{\it t\,}}(r(\rho))=\rho$. So the spanning set (\ref{1.1})
can be reduced to the spanning set
\begin{equation}\label{1.2}
r(\rho)u(\pi)v_\Lambda,\qquad \rho\in \ell \!\text{{\it t\,}}\big(\bar R\big), \ \pi\in\mathcal P.
\end{equation}

In order to reduce this spanning set to a basis, one needs
relations among vectors of the form $r(\rho)u(\pi) v_\Lambda$, or
relations among operators of the form $r(\rho)u(\pi)$, in
\cite{MP2, MP3} they are called relations among relations. Their
precise form is given by (\ref{2.19}), let us write here simply
\begin{equation}\label{1.3}
r(\rho)u(\pi)\sim r(\rho')u(\pi')\qquad\text{if}
\quad \rho\pi= \rho'\pi'.
\end{equation}

In \cite{MP2} all relations among relations are found for
$\tilde{\mathfrak g}={\mathfrak sl}(2,\mathbb
C)\,\widetilde{}$\,-modules. As a result, the spanning set
(\ref{1.2}) is reduced to a spanning subset parameterized by
colored partitions in $\mathcal P$ of the form $\rho\pi$, where
$\rho\in \ell \!\text{{\it t\,}}\big(\bar R\big)$, $\pi\in\mathcal
P$. The linear independence of such set of vectors in $M(\Lambda)$
is easy to prove.

In a sharp contrast to the general results on algebraic properties
of relations $\bar R$, the combinatorial relations (\ref{1.3})
among these relations are very difficult to handle even for
$\tilde{\mathfrak g}={\mathfrak sl}(2,\mathbb
C)\,\widetilde{}$\,-modules. In the present paper some arguments
and results in \cite[Chapters 6 and 8]{MP2} and \cite{MP3} are
formulated and extended for general affine Lie algebras. The main
result is Theorem~\ref{T:Theorem 2.8}, guaranteeing the existence
of combinatorial relations among relations (\ref{1.3}), provided
that
\begin{equation}\label{1.4}
\sum_{\pi\in\mathcal P^\ell(n)} N(\pi)=\dim Q(n).
\end{equation}
Here the left hand side is certain purely combinatorial quantity
defined for colored partitions in terms of leading terms $\ell
\!\text{{\it t\,}}(\bar R)$, and $\mathcal P^\ell(n)$ denotes the
set of colored partitions of degree $n$ and length $\ell$. On the
other side, $Q(n)$ is the space of coefficients of degree $n$ of
vertex operators $Y(u\otimes v,z)=Y(u,z)\otimes Y(v,z)$ associated
to elements of $Q$. By assumption $Q$ should be in the kernel of
the map
$$
\Phi \colon \bar R {\mathbf 1}\otimes N(k\Lambda_0)\rightarrow
N(k\Lambda_0),\quad \Phi (u\otimes v)=u_{-1}v,
$$
and such that $Q\subset (\bar R {\mathbf 1}\otimes
N(k\Lambda_0))_\ell$ and  $\ell \!\text{{\it t\,}}
\big(Q(n)\big)\subset\mathcal P^\ell(n)$. Roughly speaking, the
left hand side of (\ref{1.4}) counts the number of relations
(\ref{1.3}) needed for combinatorial arguments, and the right hand
side counts the number of  relations obtained by using
representation theory of vertex operator algebras.

The proof of Theorem~\ref{T:Theorem 2.8} follows the proof of
Lemma 4 in \cite{MP3}, the novelty is a general construction of
relations among relations by using vertex operators $Y(q,z)$
associated to elements $q\in \ker\Phi$. Although the result holds
under very mild assumptions on $R$, applications are expected
mainly for standard modules of affine Lie algebras in a manner
described above.

By Theorem~\ref{T:Theorem 2.8} the construction of combinatorial
basis of $M^1(\Lambda)$ within the spanning set (\ref{1.2}) is
reduced to verifying the equality (\ref{1.4}). It is hoped that
the equality holds in some generality, for now there are only few
examples: Besides the motivating example for ${\mathfrak
sl}(2,\mathbb C)\,\widetilde{}$\,, Ivica Siladi\' c has found in
\cite{S} a basis $\tilde B$ such that for a level 1 twisted
${\mathfrak sl}(3,\mathbb C)\,\widetilde{}$\,-module the equality
holds. For another motivating example, the basic ${\mathfrak
sl}(3,\mathbb C)\,\widetilde{}$\,-module in the homogeneous
picture, with a basis $\tilde B$ involving root vectors, the
equality (\ref{1.4}) does not hold for certain weight subspaces of
$Q$, and one is forced to perform a sort of Gr\"obner base theory
procedure to obtain a combinatorial basis of $M^1(\Lambda)$. It
might be that here Theorem~\ref{T:Theorem 2.8} could be applied as
well if one takes a somewhat larger generating subspace $R\subset
N^1(k\Lambda_0)$.

Throughout the paper $k\in\mathbb C$ is fixed, it will be the level of representations.

Some ideas worked out in this paper come {}from collaboration with
Arne Meurman for many years, and some even further {}from
collaboration with Jim Lepowsky. I thank both Arne Meurman and Jim
Lepowsky for their implicit contribution to this work. I also
thank Ivica Siladi\' c for numerous stimulating discussions.


\section{Relations among relations in terms of vertex operators and leading terms}

\subsection*{Bases of twisted affine Lie algebras and colored partitions}

Let ${\mathfrak g}$ be a simple complex Lie algebra and $\sigma$ an
automorphism of ${\mathfrak g}$ of finite order $T$, so that
$\sigma^T=1$. Set $\varepsilon=\exp(\tfrac{2\pi\sqrt{-1}}{T})$ and
${\mathfrak g}_{[j]}=\{x\in{\mathfrak g}\mid \sigma x=\varepsilon^j
x\}$ for $j\in\mathbb Z/T\mathbb Z$. Let $B_{[j]}$ be a basis of
${\mathfrak g}_{[j]}$ and
$$
B=\bigcup_{j=0}^{T-1} B_{[j]}
$$
a basis of ${\mathfrak g}$. Let $\mathfrak h$ be a Cartan
subalgebra of ${\mathfrak g}$ and $\langle \ , \ \rangle$ a
symmetric invariant bilinear form on ${\mathfrak g}$. Via this
form we identify $\mathfrak h$ and $\mathfrak h^*$ and we assume
that $\langle \theta , \theta \rangle=2$ for the maximal root
$\theta$ (with respect to some fixed basis of the root system).

Set $\hat{\mathfrak g}{[\sigma]} =\coprod_{j\in\mathbb
Z}{\mathfrak g}_{[j]}\otimes t^{j/T}+\mathbb C c$. Then
$\tilde{\mathfrak g}{[\sigma]}=\hat{\mathfrak g}{[\sigma]}+\mathbb C
d$ is the associated twisted affine Lie algebra (cf. \cite{K})
with the commutator
$$
[x(i),y(j)]=[x,y](i+j)+i\delta_{i+j,0}\langle x,y\rangle c.
$$
Here, as usual, $x(s/T)=x\otimes t^{s/T}$ for $x\in{\mathfrak
g}_{[s]}$ and $s\in\mathbb Z$, $c$ is the canonical central
element, and $[d,x(i)]=ix(i)$. Sometimes we shall denote
${\mathfrak g}_{[j]}\otimes t^{j/T}$ by ${\mathfrak g}( {j/T})$.
We define the automorphism $\sigma$ of $\tilde{\mathfrak
g}{[\sigma]}$ by $\sigma(x(s/T))=(\sigma x)(s/T)=\varepsilon^s
x(s/T)$, $\sigma c=c$, $\sigma d=d$.

In the case when the automorphism $\sigma$ is the identity on
${\mathfrak g}$, i.e., $\sigma=\text{id}$ and $T=1$, we write
$\tilde{\mathfrak g}$ instead of $\tilde{\mathfrak g}{[\text{id}]}$.
In this case we identify ${\mathfrak g}={\mathfrak g}(0)$.

\begin{remark}
As a general rule, we shall write $\mathcal X^\sigma$ for an
object whose construction involves $\sigma$, and we shall write
$\mathcal X$ instead of $\mathcal X^{\text{id}}$ when
$\sigma=\text{id}$. However, even when we consider the case
$\sigma\neq\text{id}$, we shall usually need both
$\tilde{\mathfrak g}$ and $\tilde{\mathfrak g}{[\sigma]}$ and the
corresponding objects $\mathcal X$ and $\mathcal X^\sigma$. In
that case we define the automorphism $\sigma$ of
$\tilde{\mathfrak g}$ by $\sigma(x(s))=(\sigma x)(s)$, $\sigma
c=c$, $\sigma d=d$, and the corresponding objects $\mathcal X$ and
$\mathcal X^\sigma$ will both inherit the corresponding
automorphism $\sigma$, like, for example, the enveloping algebras
$\mathcal U$ and $\mathcal U^\sigma$ bellow.
\end{remark}

Let $\mathcal U^\sigma=U(\hat{\mathfrak g}{[\sigma]})/(c-k)$,
where $U(\hat{\mathfrak g}{[\sigma]})$ is the universal enveloping
algebra of $\hat{\mathfrak g}{[\sigma]}$, and for fixed $k\in
\mathbb C$ we denoted by $(c-k)$ the ideal generated by the
element $c-k$. Note that $\tilde{\mathfrak g}{[\sigma]}$-modules
of level $k$ are $\mathcal U^\sigma$-modules. Note that
$U(\hat{\mathfrak g}{[\sigma]})$ is graded by the derivation $d$,
and so is the quotient $\mathcal U^\sigma$. The associative
algebra $\mathcal U^\sigma$ also inherits {}from $U(\hat{\mathfrak
g}{[\sigma]})$ the filtration $\mathcal U^\sigma_\ell$,
$\ell\in\mathbb Z_{\geq 0}$; let us denote by $\mathcal
S^\sigma\cong S(\bar{\mathfrak g}{[\sigma]})$ the corresponding
commutative graded algebra.

We fix a basis $\tilde B$ of $\tilde{\mathfrak g}{[\sigma]}$,
$$
\tilde B=\bar B\cup\{c, d\},\quad
\bar B=\bigcup_{j\in\mathbb Z} B_{[j]}\otimes t^{j/T},
$$
so that $\bar B$ may also be viewed as a basis of
$\bar{\mathfrak g}{[\sigma]}=\hat{\mathfrak g}{[\sigma]}/\mathbb C
c$. Let $\preceq$ be a linear order on $\bar B$ such that
$$
i<j\quad\text{implies}\quad x(i)\prec y(j).
$$

 The symmetric algebra $\mathcal S^\sigma$ has a basis $\mathcal P$
consisting of monomials in basis elements $\bar B$. Elements $\pi\in\mathcal P$ are
finite products of the form
$$
\pi=\prod_{i=1}^\ell b_i(j_i), \quad b_i(j_i)\in\bar B,
$$
and we shall say that $\pi$ is a colored partition of degree
$|\pi|=\sum_{i=1}^\ell j_i\in\frac1T \mathbb Z$ and length $\ell
(\pi)=\ell$, with parts $b_i(j_i)$ of degree $j_i$ and color
$b_i$. The set of all colored partitions of degree $n$ and length
$\ell$ is denoted as $\mathcal P^\ell(n)$. We shall usually assume
that parts of $\pi$ are indexed so that
$$
b_1(j_1)\preceq b_2(j_2)\preceq \dots\preceq b_\ell(j_\ell).
$$
The basis element $1\in\mathcal P$ we call a colored
partition of degree 0 and length 0, we may also denote it by
$\emptyset$, suggesting it has no parts. Note that $\mathcal P\subset\mathcal S^\sigma$ is a
monoid with the unit element 1, the product of monomials $\pi$ and $\rho$ is
denoted by $\pi\rho$.

We shall fix the monomial basis
\begin{equation}\label{2.1}
u(\pi)=b_1(j_1)b_2(j_2) \dots b_n(j_\ell), \quad \pi\in\mathcal P,
\end{equation}
of the enveloping algebra $\mathcal U^\sigma$. Then we have a slight modification
of Lemma 6.4.1 in \cite{MP2}:

\begin{lemma}
\label{T:Lemma 2.1} Let $\pi \in \mathcal P $. Then there exist a
restricted $\tilde{{\mathfrak g}}[\sigma]$-module $M$ of level
$k$, a vector $v_\pi \in M$ and a functional $v^* _\pi \in M^*$
such that
\begin{gather*}
\langle v^* _\pi, u(\pi)v_\pi \rangle \ne 0,\\
\langle v^* _\pi, u(\pi')v_\pi \rangle = 0 \quad \text{if}
\quad
\pi' \in \mathcal P, \pi' \ne \pi, \ell(\pi') \le
\ell(\pi).
\end{gather*}
\end{lemma}

\begin{proof}
Denote by $x^*$, $x\in B$, the elements of the dual basis of $B$ with
respect to $\langle\ ,\ \rangle$.
For an element $b =x(i)\in \bar B$ define $\bar b=x^*(-i)$. Note that
for $x(i), y(j)\in\bar B$ we have
\begin{equation}\label{2.2}
[x(i),y^*(-j)]=[x,y^*](i-j)+i\delta_{i,j}\delta_{x,y} c.
\end{equation}
Set
$$
\tilde{\mathfrak g}{[\sigma]}_{>0}
=\coprod_{j>0}{\mathfrak g}_{[j]}\otimes t^{j/T}+\mathbb C
d,\qquad\tilde{\mathfrak g}{[\sigma]}_{\geq 0} =\coprod_{j\geq
0}{\mathfrak g}_{[j]}\otimes t^{j/T}+\mathbb C d.
$$
For $s\in\mathbb C$ denote by $\mathbb C v_s$ the one-dimensional
$(\tilde{\mathfrak g}{[\sigma]}_{\geq 0}+\mathbb C c)$-module on
which $\tilde{\mathfrak g}{[\sigma]}_{\geq 0}$ acts trivially and
$c$ as multiplication by $s$. Define induced modules
$$
W(s)=U(\tilde{\mathfrak g}{[\sigma]})\otimes_{U(\tilde{\mathfrak g}{[\sigma]}_{\geq
0}+\mathbb C c)}\mathbb C v_s, \qquad
W=U(\tilde{\mathfrak g}{[\sigma]})/U(\tilde{\mathfrak g}{[\sigma]}_{>0}+\mathbb
C c).
$$
Let $\pi = \pi_- \pi_0 \pi_+$, where
$\pi_- \in \mathcal P $ has parts of negative degrees, $\pi_0\in \mathcal P $ has
parts of degree zero,  and $\pi_+\in \mathcal P $ has parts of positive degrees.
Let $p=\ell(\pi_-)$, $\ell(\pi_0)=q$, $r=\ell(\pi_+)$, $\ell=\ell(\pi)=p+q+r$
and set
$$
M=W(k-r)\otimes(\otimes^p _{i=1} W(0))\otimes W
\otimes (\otimes^r _{i=1} W(1)).
$$
Then $M$ is a $\tilde{\mathfrak g}{[\sigma]}$-module of level $k$.
Obviously $M$ is a restricted
$\tilde{\mathfrak g}{[\sigma]}$-module, i.e., for every $v\in M$
there is $n$ such that $({\mathfrak g}_{[j]}\otimes t^{j/T})v=0$
for $j\geq n$.

Recall that $\mathcal U^\sigma$ is a quotient of
$U(\hat{\mathfrak g}{[\sigma]})$ and that $u(\rho)\in\mathcal
U^\sigma$. Here we define elements $u(\rho)\in
U(\hat{\mathfrak g}{[\sigma]})$, $\rho\in\mathcal P$, by the same
formula (\ref{2.1}). In $W(0)$ we choose the basis $u(\rho) v_0$,
where $\rho \in \mathcal P$ has parts of negative degrees; for
$W(k-r)$ and $W$ we choose bases in a similar way. For $\rho$ with
parts $b_1 \preceq \cdots \preceq  b_n$ of positive degrees we set
$\bar u(\bar\rho)=\bar b_1 \dots \bar b_n$, and in $W(1)$ we
choose the basis of the form $\bar u(\bar\rho) v_1$. For a basis
of $M$ we choose the corresponding tensor products.

Let $\pi_+ = b_1\dots b_r$ and let $v_\pi\in M$ be the basis vector
$$
v_\pi=v_{k-r}\otimes(\otimes^p _{i=1} v_0)\otimes 1 \otimes
(\otimes^r _{i=1} \bar b_i v_{1}).
$$
Let $v^* _\pi$ be the element of dual basis of $M^*$ corresponding to
the vector
$$
v=v_{k-r}\otimes(\otimes^p _{i=1} d_i v_0) \otimes u(\pi_0)1 \otimes
(\otimes^r _{i=1} v_{1}),
$$
where $d_i \in \bar B$ are the parts of $\pi_- = d_1\dots d_p$. It
is clear {}from our choice and (\ref{2.2}) that
$$
v^* _\pi(u(\pi)v_\pi) > 0.
$$

Now let $\pi' \in \mathcal P (\bar B)$, $\ell (\pi') \le \ell$,
and
assume that $u(\pi')v_\pi$ has a nontrivial component along
the basis vector
$v$. Write $\pi'=\pi' _- \pi' _0 \pi' _+$, notation being as
before. In
order to have a component $C v = C(v_{k-r}\otimes
(\otimes^p _{i=1} d_i v_0) \otimes u(\pi_0)1 \otimes
(\otimes^r _{i=1} v_{1}))$, $C \ne 0$, $\pi' _-$ must
contain all parts $d_i$ of $\pi$ and
$\pi' _0$ must contain all parts of $\pi _0$. Hence $\ell(\pi' _+)\le r$. But to
have  $u(\pi' _+)(\otimes^r _{i=1} \bar b_i v_{1})=
C_1(\otimes^r _{i=1}  v_{1})+ \cdots$, $C_1 \ne 0$, we need
$\ell(\pi' _+) \ge r$. Hence $\ell(\pi' _-)= r$, $\ell(\pi' _0)=q$,
$\ell(\pi' _+) = r$. But then $\pi' _- = \pi_-$, $\pi' _0= \pi_0$ and,
because of (\ref{2.2}), $\pi'_+ = \pi_+$. Hence  $\pi' = \pi$.
\end{proof}

\subsection*{A filtration on a completion of the enveloping algebra
and leading terms}

The coefficients of vertex operators may be seen as infinite sums
of elements in $\mathcal U^\sigma$ and in certain arguments it is
necessary to use some completion $\overline{\mathcal U^\sigma}$
where these infinite sums converge. In \cite{MP1, MP2} the
topology of pointwise convergence on a certain set (category) of
$\mathcal U$-modules was used. However, it seems that another kind
of topology, used in \cite{H}, \cite{FF}, \cite{FZ} and \cite{KL},
for example, has several advantages. Here we shall adopt the point
of view taken by I.~B.~Frenkel and Y.~Zhu: Let us denote the
homogeneous components of the graded algebra $\mathcal U^\sigma$
by $\mathcal U^\sigma(n)$, $n\in\frac1T \mathbb Z$. We take
$$
V_p(n)=\sum_{i\geq p} \mathcal U^\sigma(n-i)\mathcal U^\sigma(i),\qquad p=1,2,\dots \ ,
$$
to be a fundamental system of neighborhoods of $0\in \mathcal U^\sigma(n)$.
It is easy to see that
$$
\bigcap_{p=1}^\infty V_p(n)=\{0\},
$$
so we have a Hausdorff topological group $(\mathcal U^\sigma(n),+)$, and we
denote by $\overline{\mathcal U^\sigma(n)}$ the corresponding completion.
Then
$$
\overline{\mathcal U^\sigma}=\coprod_{n\in\frac1T \mathbb Z}\overline{\mathcal U^\sigma(n)}
$$
is a topological ring. It is clear that we have the action of $\sigma$ on both
$\overline{\mathcal U}$ and $\overline{\mathcal U^\sigma}$.

Clearly $\bar B\subset\mathcal P$, viewed as
colored partitions of length 1. We assume that on $\mathcal P$ we have
a linear order $\preceq$ which extends the order $\preceq$ on
$\bar B$. Moreover, we assume that order $\preceq$ on $\mathcal P$ has the following
properties:
\begin{itemize}
\item  $\ell(\pi)>\ell(\kappa)$ implies $\pi\prec\kappa$.
\item  $\ell(\pi)=\ell(\kappa)$, $|\pi|<|\kappa|$ implies $\pi\prec\kappa$.
\item  Let $\ell(\pi)=\ell(\kappa)$, $|\pi|=|\kappa|$. Let $\pi$ be a partition
$b_1(j_1)\preceq b_2(j_2)\preceq \dots\preceq b_\ell(j_\ell)$ and
$\kappa$ a partition
$a_1(i_1)\preceq a_2(i_2)\preceq \dots\preceq a_\ell(i_\ell)$.
Then $\pi\preceq \kappa$ implies $j_\ell\leq i_\ell$.
\item  Let $\ell \ge 0$, $n \in \frac1T \mathbb Z$ and let
$S \subset \mathcal P$ be a nonempty subset such that
all $\pi$ in $S$
have length $\ell(\pi) \le \ell$ and degree $\vert \pi \vert
= n$. Then $S$
has a minimal element.
\item   $\mu  \preceq  \nu$ implies
$\pi\mu  \preceq \pi\nu$.
\end{itemize}

\begin{remark}
An order with these properties is used in \cite{MP2, MP3}; the
first three properties are a part of definition, the last two
properties are guaranteed by Lemmas 6.2.1 and 6.2.2 in \cite{MP2}.
\end{remark}

For $\pi\in \mathcal P$, $|\pi|=n$, set
$$
U_{[\pi]}^\mathcal P =\overline{\mathbb C\text{-span}\{u(\pi')\mid |\pi'|=|\pi|, \pi' \succeq\pi\}}
\subset\overline{\mathcal U^\sigma(n)},
$$
$$
U_{(\pi)}^\mathcal P  =\overline{\mathbb C\text{-span}\{u(\pi')\mid |\pi'|=|\pi|, \pi' \succ \pi\}}
\subset\overline{\mathcal U^\sigma(n)},
$$
the closure taken in $\overline{\mathcal U^\sigma(n)}$. Set
$$
U^\mathcal P (n)=  \bigcup_{\pi\in \mathcal P, \ |\pi|=n}
U_{[\pi]}^\mathcal P ,\qquad
U^\mathcal P =\coprod_{n\in\frac1T\mathbb Z}U^\mathcal P (n) \subset\overline{\mathcal U^\sigma}.
$$
As suggested by our notation, the construction of $U^\mathcal P$
depends on our choice of $(\mathcal P, \preceq)$. Since by
assumption $\mu  \preceq  \nu$ implies $\pi\mu  \preceq \pi\nu$,
we have that $U^\mathcal P $ is a subalgebra of
$\overline{\mathcal U^\sigma}$. Obviously $U^\mathcal P (n)$ is
filtered by $U_{[\pi]}^\mathcal P $. As in \cite[Lemma
6.4.2]{MP2}, by using Lemma~\ref{T:Lemma 2.1} we obtain:

\begin{lemma}
\label{T:Lemma 2.2} For $\pi \in \mathcal P $ we have
$U_{[\pi]}^\mathcal P =\mathbb C u(\pi)+ U_{(\pi)}^\mathcal P $.
Moreover,
$$
\dim U_{[\pi]}^\mathcal P / U_{(\pi)}^\mathcal P = 1.
$$
\end{lemma}
\noindent
For $u \in U_{[\pi]}^\mathcal P $, $u \notin U_{(\pi)}^\mathcal P $ we define the
{\it leading term}
$$
\ell \!\text{{\it t\,}}(u) = \pi.
$$

\begin{proposition}
\label{T:Proposition 2.3}  Every element
 $u \in U^\mathcal P (n)$, $u \neq 0$, has a unique leading
 term $\ell \!\text{{\it t\,}}(u)$.
\end{proposition}

\begin{proof}
Since in the construction of $\overline{\mathcal U^\sigma(n)}$ we
used a countable fundamental system of neighborhoods of $0$, we
may work only with Cauchy sequences, and step (3) in the proof of
\cite[Proposition 6.4.3]{MP2} can be almost literally applied: The
main point to prove is that
$$
\bigcap_{\ell(\pi)=\ell, \ |\pi|=n} U_{[\pi]}^\mathcal P  =
\bigcup_{\ell(\pi) < \ell, \ |\pi|=n} U_{[\pi]}^\mathcal P .
$$
The inclusion $\supset$ is obvious {}from the assumed properties
of our order $\preceq$. So let $x \in U_{[\pi]}^\mathcal P $ for
all $\pi$ such that $\ell (\pi) = \ell$, $|\pi|=n$. For each $\pi$
we may choose a sequence $(x_i ^\pi)$ in  $\mathbb
C\text{-span}\{u(\pi')\mid |\pi'|=n, \pi' \succeq\pi\}$ such that
$x_i ^\pi \rightarrow x$. Let
$$
x_i ^\pi  = \sum_{\ell(\kappa)=\ell, \ |\kappa|=n} C^\pi _{\kappa ,i}
u(\kappa) +
x^\pi _{i,\ell - 1}
= x^\pi _{i,\ell} + x^\pi _{i,\ell - 1},
$$
where $x^\pi _{i,\ell - 1} \in  \mathcal U^\sigma_{\ell- 1}(n)$. For a
sequence $\pi^{(1)} \prec \pi^{(2)} \prec \cdots$ such that
$\ell(\pi^{(s)})=\ell$, $|\pi^{(s)}|=n$
we can construct a diagonal sequence
$$
x^{\pi^{(s)}} _{i_s}\in x+V_s(n), \ s=1,2,\dots ,\quad
\quad i_s>i_{s'} \ \text{for} \ s>s',
$$
so that $x^{\pi^{(s)}} _{i_s}\rightarrow x$ when $s\rightarrow \infty$.

Note that for a given $M$ there is only finitely many partitions
$n=j_1+\dots +j_\ell$, $j_p\in\frac1T\mathbb Z$, $j_1\leq\dots \leq j_\ell$,
such that $j_\ell\leq M$, and these partitions can be colored with elements
in $B$ in only finitely many ways. Hence for
$$
\pi^{(s)}=\bigl(b^{(s)}_1(j^{(s)}_1)\preceq b^{(s)}_2(j^{(s)}_2)
\preceq \dots\preceq b^{(s)}_\ell (j^{(s)}_\ell)\bigr)
$$
we have that $j^{(s)}_\ell\rightarrow \infty$ when $s\rightarrow
\infty$. Now it follows {}from our construction and the assumed
properties of our order $\preceq$ that $x^{\pi^{(s)}} _{i_s,\ell}
\rightarrow 0$ when $s\rightarrow \infty$. Hence we have the
sequence
$$
x_s = x_{i_{s},\ell -1} ^{\pi^{(s)}} \in U_{[\pi_0]}^\mathcal P \subset
\bigcup_{\ell(\pi)<\ell} U_{[\pi]}^\mathcal P
$$
such that $x_s \rightarrow x$, where we denoted by $\pi_0$ the
minimal element among colored partitions $\pi$
of length $\ell(\pi) \le \ell-1$ and degree $\vert \pi \vert= n$.
So the other inclusion $\subset$  holds as well.

Let $x \in U_{[\pi]}^\mathcal P $. Assume that $x \notin U_{[\pi']}^\mathcal P $,
where $\pi' \succ \pi$,
$\ell(\pi')=\ell(\pi)=\ell$, $|\pi'|=|\pi|=n$. Arguing as before
we see that the interval $[\pi,\pi']$ is finite, so there is
$\pi \preceq \tau \preceq \pi'$ such that
$x \in U_{[\tau]}^\mathcal P  \backslash U_{(\tau)}^\mathcal P $, and $x$ has a
leading term $\tau$. On the other hand, if $x \in U_{[\pi']}^\mathcal P $
for all $\pi' \succ \pi$,
$\ell(\pi')=\ell$, $|\pi'|=n$, then we apply the first part of the argument,
and in a finite number of steps we see that either $x=0$ or
$x$ has a leading term.

Since for $\pi\prec\kappa$ we have
$U_{(\pi)}^\mathcal P \supset U_{[\kappa]}^\mathcal P $, it is
impossible to have  $u \in U_{[\pi]}^\mathcal P $, $u \notin U_{(\pi)}^\mathcal P $ and
$u \in U_{[\kappa]}^\mathcal P $, $u \notin U_{(\kappa)}^\mathcal P $.
Hence the leading term is unique.
\end{proof}

By Proposition~\ref{T:Proposition 2.3} every nonzero homogeneous
$u$ has the unique leading term. For a nonzero element $u\in
U^\mathcal P $ we define the leading term $\ell \!\text{{\it
t\,}}(u)$ as the leading term of the nonzero homogeneous component
of $u$ of smallest degree. For a subset $S\subset U^\mathcal P $
set
$$
\ell \!\text{{\it t\,}}(S)=\{\ell \!\text{{\it t\,}}(u)\mid u\in S, u\neq 0\}.
$$
Clearly Lemma~\ref{T:Lemma 2.2} and Proposition~\ref{T:Proposition
2.3} imply the following:
\begin{proposition}
\label{T:Proposition 2.4} For all $u, v \in U^\mathcal P
\backslash\{0\}$ we have $\ell \!\text{{\it t\,}}(uv)=\ell
\!\text{{\it t\,}}(u)\ell \!\text{{\it t\,}}(v)$.
\end{proposition}
\begin{proposition}
\label{T:Proposition 2.5} Let $W\subset U^\mathcal P $ be a finite
dimensional subspace and let \newline
$\ell \!\text{{\it
t\,}}(W)\rightarrow W$ be a map such that
$$
\rho\mapsto w(\rho),\qquad \ell \!\text{{\it t\,}}\big( w(\rho)\big)=\rho.
$$
Then $\{w(\rho)\mid \rho\in\ell \!\text{{\it t\,}}(W)\}$ is a basis of $W$.
\end{proposition}

\subsection*{Vertex operators with coefficients in a completion
 of the enveloping algebra}

The untwisted affine Lie algebra $\tilde{\mathfrak g}$ gives rise
to the universal vertex operator algebra
$N(k\Lambda_0)=U(\tilde{\mathfrak g})\otimes_{U(\tilde{\mathfrak
g}_{\geq 0}+\mathbb C c)}\mathbb C v_k$ for $k\neq -g^\vee$, where
$g^\vee$ is the dual Coxeter number of ${\mathfrak g}$ (see
\cite{FZ} and \cite{L1}, we use the notation {}from \cite{MP2});
it is generated by fields
\begin{equation}\label{2.3}
x(z)=\sum_{n\in\mathbb Z}x_nz^{-n-1},\qquad x\in{\mathfrak g},
\end{equation}
where we set $x_n=x(n)$ for $x\in{\mathfrak g}$. As usual, we
shall write $Y(v,z)=\sum_{n\in\mathbb Z}v_nz^{-n-1}$ for the
vertex operator (field) associated with a vector $v\in
N(k\Lambda_0)$. By formal definition, the coefficients $v_n$ are
linear operators on $N(k\Lambda_0)$, but sometimes, as in
\cite{FF}, \cite{FZ} and \cite{MP2, MP3}, it is convenient to
think of $v_n$ as elements in a completion of the universal
enveloping algebra of $\tilde{\mathfrak g}$, i.e.,
$v_n\in\overline{\mathcal U}$. Then for any restricted
$\tilde{\mathfrak g}$-module $M$ the elements
$v_n\in\overline{\mathcal U}$ act on $M$. The construction of a
map
\begin{equation*}
Y \colon N(k\Lambda_0) \rightarrow \overline{\mathcal
U}[[z,z^{-1}]],\qquad v\mapsto Y(v,z)=\sum_{n\in\mathbb
Z}v_nz^{-n-1},
\end{equation*}
was given by I.~B.~Frenkel and Y.~Zhu in \cite[Definition
2.2.2]{FZ}. Their result can be interpreted in another way by
using the approach developed by Haisheng Li: As in \cite[Section
2.2]{FZ}, we may consider regular and mutually local fields $a(z)$
and $b(z)$ with coefficients $a_n, b_n\in\overline{\mathcal U}$
and, as in \cite[Lemma 3.1.4]{L1}, we may define products $a(z)_m
b(z)$ for $m\in\mathbb Z$. These products are well defined because
infinite sums, appearing as coefficients of $a(z)_m b(z)$, are
convergent in $\overline{\mathcal U}$ (cf. \cite[Definition
2.2.1]{FZ}). Since the proofs of \cite[Propositions 3.2.7 and
3.2.9]{L1} apply literally, we may introduce the notion of local
system of vertex operators with coefficients in
$\overline{\mathcal U}$ for which Theorem 3.2.10 in \cite{L1}
applies. In particular, we have a vertex operator algebra $V$
generated by fields (\ref{2.3}). Again by literally repeating the
arguments in \cite{L1}, we get that Frenkel-Zhu's map
$N(k\Lambda_0)\rightarrow V$ is an isomorphism of vertex operator
algebras. To be more precise, by \cite[Section 2.2]{FZ} $V$ is a
highest weight $\tilde {\mathfrak g}$-module, and we have a
surjective homomorphism $N(k\Lambda_0)\rightarrow V$ of vertex
operator algebras; it must be an isomorphism because the adjoint
representation of vertex operator algebra is faithful.

In the case when $\sigma$ is not the identity on ${\mathfrak g}$,
we may use Haisheng Li's theory \cite{L2} of local systems of
twisted vertex operators. The twisted vertex operators
\begin{equation*}
x(z)=\sum_{n\in\frac1T\mathbb Z}x_nz^{-n-1},\qquad
x\in{\mathfrak g},
\end{equation*}
where we set $x_{s/T}=x(s/T)$ for $x\in{\mathfrak g}_{[s]}$,
generate a vertex operator algebra $V$, the $m$-products of fields
are given by \cite[Definition 3.7]{L2}. The results in \cite{L2}
also imply that we have a surjective homomorphism
$N(k\Lambda_0)\rightarrow V$ of vertex operator algebras given by
a map
\begin{equation}\label{2.6}
Y^\sigma \colon N(k\Lambda_0)\rightarrow \overline{\mathcal
U^\sigma}[[z^{\frac1T},z^{-\frac1T}]],\qquad v\mapsto
Y^\sigma(v,z)=\sum_{n\in {\frac1T}\mathbb Z}v_nz^{-n-1}.
\end{equation}
I am not aware of the existence of a faithful $\sigma$-twisted
representation of $N(k\Lambda_0)$ in general (see \cite[Section
5]{L2} for inner automorphisms) and I don't know whether
(\ref{2.6}) is an isomorphism in general. It is clear that for any
restricted $\tilde{\mathfrak g}{[\sigma]}$-module $M$ the action
of coefficients of $Y^\sigma(v,z)$ on $M$ gives rise to the
twisted representation $Y^\sigma_M(v,z)$ of $N(k\Lambda_0)$ on
$M$.

When $\sigma\neq \text{id}$, we have the action of automorphism $\sigma$
on the vertex operator algebra $N(k\Lambda_0)$, as well as on both
$\overline{\mathcal U}$ and $\overline{\mathcal U^\sigma}$. It is easy to see that
for both untwisted and twisted fields we have $(\sigma v)_n=\sigma (v_n)$.

By following the notation in \cite{FF} we set
$$
U_{\text{loc}}^\sigma=\mathbb C\text{-span}\{v_n \mid
v \in N(k\Lambda_0), n \in \tfrac1T\mathbb Z\}\subset\overline{\mathcal U^\sigma},
$$
where $v_n$ denotes a coefficient in $Y^\sigma(v,z)$. {}From the
commutator formula (cf. \cite[(2.31)]{L2}) we see that
$U_{\text{loc}}^\sigma$ is a Lie algebra. Let us denote by
$U^\sigma$ the associative subalgebra of $\overline{\mathcal
U^\sigma}$ generated by $U_{\text{loc}}$. Clearly
$$
U^\sigma=\coprod_{n\in\frac1T\mathbb Z} U^\sigma(n),
$$
where $U^\sigma(n)\subset U^\sigma$ is the homogeneous subspace of degree $n$.

{}From the normal order product formula (cf. \cite[(3.6)]{L2}) we
see by induction that
$$
U_{\text{loc}}^\sigma \subset
U^\mathcal P =\coprod_{n\in\frac1T\mathbb Z}\Big( \bigcup_{
\substack{ |\pi|=n}} U_{[\pi]}^\mathcal P \Big).
$$
Since $U^\mathcal P\subset \overline{\mathcal U^\sigma}$ is a
subalgebra, we have $U^\sigma \subset U^\mathcal P$. Hence
Propositions~\ref{T:Proposition 2.3} and \ref{T:Proposition 2.4}
imply that every nonzero element in $U^\sigma$ has the leading
term and that $\ell \!\text{{\it t\,}}(uv)=\ell \!\text{{\it
t\,}}(u)\ell \!\text{{\it t\,}}(v)$. Note that $\ell(\ell
\!\text{{\it t\,}}(uv))>\ell(\ell \!\text{{\it t\,}}([u,v]))$, so
in particular $\ell \!\text{{\it t\,}}(uv)\prec\ell \!\text{{\it
t\,}}([u,v])$.

We may summarize our constructions with the following:
\begin{proposition}
\label{T:Proposition 2.6}
$$
\mathcal U^\sigma\subset U^\sigma \subset U^\mathcal P \subset \overline{\mathcal U^\sigma}.
$$
\end{proposition}

\subsection*{A combinatorial formulation of relations among relations}

{}From now on we fix $k\neq -g^\vee$. Then $N(k\Lambda_0)$ is a
vertex operator algebra, and we denote by $\omega$ the conformal
vector and by $L_n$, $n\in\mathbb Z$, the elements of the Virasoro
algebra. We assume that $N(k\Lambda_0)$ is reducible and we denote
by $N^1(k\Lambda_0)$ the maximal ideal.

{}From now on we also assume that $R\subset N(k\Lambda_0)$ is a
nontrivial subspace such that:
\begin{itemize}
\item  $R$ is finite dimensional,
\item  $R$ is invariant for $L_0$ and $\sigma$,
\item  $R$ is invariant for $\coprod_{n\in\mathbb Z_{\geq 0}}{\mathfrak g}(n)$,
\item  $R\subset N^1(k\Lambda_0)$.
\end{itemize}
Set
$$
\bar R^\sigma = \mathbb C\text{-span}\{r_n \mid
r \in R, n \in \tfrac1T\mathbb Z\}\subset\overline{\mathcal U^\sigma},
$$
where $r_n$ denotes a coefficient in $Y^\sigma(r,z)$. Then $\bar
R^\sigma$ is a $\tilde{{\mathfrak g}}[\sigma]$-module for the
adjoint action given by the commutator formula
$$
[x_m,r_n]=\sum_{i\geq 0}\binom{m}{i}(x_i r)_{m+n-i},\qquad
x\in{\mathfrak g}, \ r\in R.
$$
We shall say that $\bar R^\sigma$ is a loop module; in general it is reducible.
When $\sigma$ is the identity, we shall write $\bar R$ instead of
${\bar R}\sp{\text{id}}$.

Let $M(\Lambda)$ be a Verma module for
$\tilde{{\mathfrak g}}[\sigma]$. Set
$$
W(\Lambda)=\bar R^\sigma M(\Lambda).
$$
Then $W(\Lambda)$ is a $\tilde{{\mathfrak g}}[\sigma]$-submodule of
$M(\Lambda)$. If  $W(\Lambda)\neq M(\Lambda)$, then we have a
highest weight $\tilde{{\mathfrak g}}[\sigma]$-module
$$
M(\Lambda)/W(\Lambda)
$$
which is, by construction, annihilated by all $Y^\sigma(v,z)$ for
elements $v$ in the ideal ${\bar R}N(k\Lambda_0)$. We shall call
such vertex operators $Y^\sigma(v,z)$ the {\it annihilating
fields} of the $\tilde{{\mathfrak g}}[\sigma]$-module
$M(\Lambda)/W(\Lambda)$. Although obvious, it should be emphasized
that it is the coefficients of annihilating fields that annihilate
the module, and, since they define the module, we shall call them
{\it relations} for $M(\Lambda)/W(\Lambda)$.

\begin{remark}
In the case of standard modules, i.e., when $k$ is a positive
integer, we take a ${\mathfrak g}$-module $R$ generated by the
singular vector $x_\theta(-1)^{k+1}\mathbf 1$. Then ${\bar
R}N(k\Lambda_0)$ is the maximal ideal by the result of V.~G.~Kac
\cite[Corollary 10.4]{K}. Moreover, for $\Lambda$ integral
dominant $M(\Lambda)/W(\Lambda)$ is the standard module
$L(\Lambda)$, and the structure of relations for $L(\Lambda)$ is
well understood (cf. \cite{DL}, \cite{FZ}, \cite{L1} and
\cite[Chapter 5]{MP2}). This is also true in a twisted case (cf.
\cite[Theorem 4.11]{B} and \cite[Section 5]{L2}).

Our discussion might be also applicable to other admissible
representations \cite{KW}, like the series
$k\in\tfrac12 +\mathbb N$ for affine Lie algebras of type $C_n^{(1)}$
studied in \cite{A},
where all irreducible $L(k\Lambda_0)$-modules are of
the form $L(\Lambda)=M(\Lambda)/W(\Lambda)$ for admissible weights $\Lambda$ of
level $k$.
\end{remark}

Since by assumption $R$ is finite dimensional, the space $\bar
R^\sigma\subset U^\sigma$ is a direct sum of finite dimensional
homogeneous subspaces. Hence Proposition~\ref{T:Proposition 2.5}
implies that we can parametrize a basis of $\bar R^\sigma$ by the
set of leading terms $\ell \!\text{{\it t\,}} (\bar R^\sigma)$: we
fix a map
$$
\ell \!\text{{\it t\,}} (\bar{R}^\sigma) \rightarrow \bar R^\sigma ,\quad\rho \mapsto r(\rho)\quad
\text{such that}\quad r(\rho)\in U^\sigma(\vert \rho \vert),\ \ell \!\text{{\it t\,}} (r(\rho)) = \rho,
$$
then $\{r(\rho)\mid \rho\in \ell \!\text{{\it t\,}} (\bar R^\sigma)\}$ is a basis of $\bar R^\sigma$.
We will assume that this map is such that the coefficient $C$ of
``the leading term'' $u(\rho)$ in ``the expansion'' of
$r(\rho)=Cu(\rho)+\dots$ is chosen to be $C=1$. Note that our assumption
$R\subset N^1(k\Lambda_0)$ implies that $1\not\in \ell \!\text{{\it t\,}} (\bar R^\sigma)$
and that $ \ell \!\text{{\it t\,}} (\bar R^\sigma)\cdot\mathcal P$ is a proper ideal in
the semigroup $\mathcal P$.

Since the character of Verma module $M(\Lambda)$ is easily
described, a combinatorial description of the character of
$M(\Lambda)/W(\Lambda)$ may be obtained {}from a combinatorial
description of $W(\Lambda)$. So in Lepowsky-Wilson's approach a
construction of combinatorial basis of $M(\Lambda)/W(\Lambda)$ can
be obtained by constructing first a combinatorial basis of
$W(\Lambda)$. If $v_\Lambda\in M(\Lambda)$ is a highest weight
vector, then it is natural to seek for a combinatorial basis of
$W(\Lambda)$ within the spanning set
\begin{equation}\label{2.7}
r(\rho)u(\pi)v_\Lambda,\quad \rho\in \ell \!\text{{\it t\,}}(\bar R^\sigma), \pi\in\mathcal P.
\end{equation}
In order to reduce this spanning set to a basis, one needs relations
among vectors of the form $r(\rho)u(\pi)v_\Lambda$,
or relations among operators of the form $r(\rho)u(\pi)$.

For colored partitions $\kappa$, $\rho$ and $\pi=\kappa\rho$ we shall write
$\kappa=\pi/\rho$ and $\rho\subset\pi$. We shall say that $\rho\subset\pi$ is an
embedding (of $\rho$ in $\pi$), notation suggesting that $\pi$ ``contains'' all the
parts of $\rho$.
For an embedding $\rho\subset\pi$, where $\rho\in\ell \!\text{{\it t\,}} (\bar R)$, we define the
element $u(\rho\subset\pi)$ in $ U^\sigma$ by
$$
u(\rho\subset\pi)=\begin{cases} u(\pi/\rho)r(\rho) &\text{if \ } |\rho |>|\pi/\rho |,\\
r(\rho)u(\pi/\rho) &\text{if \ } |\rho |\leq |\pi/\rho
|.\end{cases}
$$

Instead of (\ref{2.7}), it is convenient to consider a slightly
modified spanning set
\begin{equation}\label{2.8}
u(\rho\subset\pi)v_\Lambda,\quad \rho\in \ell \!\text{{\it t\,}}(\bar R^\sigma), \pi\in\mathcal P.
\end{equation}
Note that, by Proposition~\ref{T:Proposition 2.4}, we have $\ell
\!\text{{\it t\,}}\big(u(\rho\subset\pi)\big)=\pi$. It is clear
that in general we shall have several embeddings in given $\pi$,
say $\rho_1 \subset \pi$, $\rho_2 \subset \pi$, \dots , giving
several vectors in the spanning set (\ref{2.8}) with the same
leading term $\ell \!\text{{\it
t\,}}\big(u(\rho_1\subset\pi)\big)=\ell \!\text{{\it
t\,}}\big(u(\rho_2\subset\pi)\big)=\dots=\pi$. If we can prove
that for any two embeddings $\rho_1 \subset \pi$, $\rho_2 \subset
\pi$ there is a relation of the form
\begin{equation}\label{2.9}
u(\rho_1 \subset \pi) \in  u(\rho_2 \subset \pi) + \overline
{\mathbb C\text{-span}\{u(\rho \subset \pi')\mid
 \rho \in \ell \!\text{{\it t\,}} (\bar{R}^\sigma), \rho \subset \pi', \pi \prec \pi'\}},
\end{equation}
then we can reduce the spanning set (\ref{2.8}) to a spanning set
of $W(\Lambda)$ indexed by colored partitions:
\begin{equation}\label{2.10}
u(\rho\subset\pi)v_\Lambda,\qquad
\pi\in\ell \!\text{{\it t\,}}(\bar R^\sigma)\cdot\mathcal P,
\end{equation}
where for each $\pi\in\ell \!\text{{\it t\,}}(\bar
R^\sigma)\cdot\mathcal P$ we choose only one embedding
$\rho\subset\pi$. Roughly speaking (i.e., if we forget the
so-called initial conditions, see \cite[Chapter 6]{MP2} for the
correct formulation), the spanning set (\ref{2.10}) should be a
basis of $W(\Lambda)$; the linear independence should follow
easily {}from the fact that our vectors
$u(\rho\subset\pi)v_\Lambda=u(\pi)v_\Lambda+\dots$ are elements of
the Verma module $M(\Lambda)$ and that the leading term $\pi=\ell
\!\text{{\it t\,}}\big(u(\rho\subset\pi)\big)$ appears only once
in the spanning set (\ref{2.10}).

We may call a relation of the form (\ref{2.9}) a combinatorial {\it relation among relations}.

\subsection*{A vertex operator algebra formulation of relations among relations}

Recall that $N(k\Lambda_0)\otimes N(k\Lambda_0)$ is a vertex
operator algebra; the fields are defined by $Y(a\otimes
b,z)=Y(a,z)\otimes Y(b,z)$, the conformal vector is $\omega\otimes
\mathbf 1+\mathbf 1\otimes \omega$ (cf.  \cite{FHL}). In
particular, the derivation $D=L_{-1}$ is given by $D\otimes
1+1\otimes D$, the degree operator $-d=L_0$ is given by
$L_0\otimes 1+1\otimes L_0$, and we have the action of ${\mathfrak
g}={\mathfrak g}(0)$ given by $x\otimes 1+1\otimes x$. We consider
the automorphism $\sigma\otimes\sigma$ of vertex operator algebra
$N(k\Lambda_0)\otimes N(k\Lambda_0)$, also denoted by $\sigma$.

For a $\sigma$-twisted $N(k\Lambda_0)$-module $M$ given by twisted
vertex operators $Y_M^\sigma(a,z)$ we have the $\sigma$-twisted
$\big(N(k\Lambda_0)\otimes N(k\Lambda_0)\big)$-module $M\otimes M$
given by $\sigma$-twisted vertex operators $Y_M^\sigma(a,z)\otimes
Y_M^\sigma(b,z)$. The coefficients $(a\otimes b)_n$ of
$Y_M^\sigma(a,z)\otimes Y_M^\sigma(b,z)$ are well defined
operators on $M\otimes M$ given by
\begin{equation}\label{2.11}
(a\otimes b)_n=\sum_{i+j+1=n}a_i\otimes b_j.
\end{equation}
Note that we have a linear map between fields
\begin{equation}\label{2.12}
Y_M^\sigma(a,z)\otimes Y_M^\sigma(b,z)
\mapsto \big(Y_M^\sigma(a,z)\big)_{-1} \big(Y_M^\sigma(b,z)\big)=Y_M^\sigma(a_{-1}b,z),
\end{equation}
where the $(-1)$-product for mutually local twisted fields is
given by \cite[Definition 3.7]{L2}. Set
\begin{equation*}
\Phi \colon N(k\Lambda_0)\otimes N(k\Lambda_0)\rightarrow
N(k\Lambda_0),\qquad \Phi(a\otimes b)=a_{-1}b.
\end{equation*}
Although we don't have a state-field correspondence for vertex
operators on the module $M$, the map $\Phi$ between vectors
corresponds to the map (\ref{2.12}) between fields. It is easy to
see that $\Phi$ intertwines the actions of $L_{-1}$, $L_{0}$,
${\mathfrak g}(0)$ and $\sigma$, so that $\ker \Phi$ is also
invariant for the actions of these operators.

The basic idea for constructing relations for annihilating fields
$Y_M^\sigma(a,z)$, $a\in\bar R N(k\Lambda_0)$ is to use the
following observation:
\begin{proposition}
\label{T:Proposition 2.6,4} \ $\sum a\otimes b\in \ker \Phi$\
implies \ $\sum \big(Y_M^\sigma(a,z)\big)_{-1}
\big(Y_M^\sigma(b,z)\big)=0$.
\end{proposition}

Since it is the coefficients of annihilating fields that are
``true'' relations, we shall need a slight refinement of the map
(\ref{2.12}). As before, we want to consider vertex operators
$Y^\sigma(a,z)$ with coefficients in $U^\sigma$. But then we have
to make sense of  (\ref{2.11}): For a fixed $n\in\frac1T\mathbb Z$
we have a linear map
$$
\chi_n^\sigma \colon N(k\Lambda_0)\otimes N(k\Lambda_0)\rightarrow
\big(U^\sigma\otimes U^\sigma)^{\tfrac1T\mathbb Z}, \qquad
a\otimes b\mapsto \left(a_i\otimes b_{n-i-1}\mid
i\in\tfrac1T\mathbb Z\right).
$$
So formally $\chi_n^\sigma(a\otimes b)$ is a sequence in
$U^\sigma\otimes U^\sigma$, but we should think of it as ``the
$n$-th coefficient $(a\otimes b)_n$ of the vertex operator
$Y^\sigma(a\otimes b,z)=Y^\sigma(a,z)\otimes Y^\sigma(b,z)$'', and
we shall formally write
$$
\chi_n^\sigma(a\otimes b)=\sum_{i+j+1=n}a_i\otimes b_j\in K_n^\sigma,
$$
where $K_n^\sigma=\chi_n^\sigma \big(N(k\Lambda_0)\otimes
N(k\Lambda_0)\big)$ is a linear subspace of $\big(U^\sigma\otimes
U^\sigma)^{\tfrac1T\mathbb Z}$. Now we define a linear map
$\Phi_n^\sigma \colon K_n^\sigma\rightarrow U^\sigma$ given by
\begin{align*}
\Phi_n^\sigma \colon \left(a_p\otimes b_{n-p-1}\mid
p\in\tfrac1T\mathbb Z\right)\mapsto &\sum_{\substack{ j\in\mathbb
Z_{<0}\\i+j+1=m}} a_{j+\tfrac rT}b_{i+\tfrac sT}+
\sum_{\substack{ j\in\mathbb Z_{\geq 0}\\i+j+1=m }} b_{i+\tfrac sT}a_{j+\tfrac rT}\\
&+\sum_{\substack{ l\in\mathbb Z_{> 0}\\ \text{finite sum}}}
(-1)^l\binom{\tfrac rT}{l} \sum_{\substack{ j=-l \\i+j+1=m }}^{-1}
\binom{l-1}{j+l} [a_{j+\tfrac rT},b_{i+\tfrac sT}]
\end{align*}
for $a\in N(k\Lambda_0)^r$, $b\in N(k\Lambda_0)^s$, $n=m+\tfrac
rT+\tfrac sT$. We should see that this map is well defined. So
first recall that we have the automorphism $\sigma$ on
$N(k\Lambda_0)$ and $U^\sigma$ and that $(\sigma
a)_p=\sigma(a_p)$. We also have actions of $\sigma\otimes 1$ and
$1\otimes\sigma$ on $N(k\Lambda_0)\otimes N(k\Lambda_0)$ and on
$\big(U^\sigma\otimes U^\sigma)^{\tfrac1T\mathbb Z}$ and
$\chi_n^\sigma$ intertwines these actions. Hence $K_n^\sigma$ is
invariant for these actions as well, and the assumption $a\in
N(k\Lambda_0)^r$, $b\in N(k\Lambda_0)^s$ means that we are
defining a map for elements in $K_n^\sigma$ with eigenvalues
$\varepsilon^r$ for $\sigma\otimes 1$ and $\varepsilon^s$ for
$1\otimes\sigma$. Next note that each term that appears on the
right hand side is a linear function of some coordinate
$a_p\otimes b_q$, $p+q+1=n$. Since infinite sums are convergent,
we have that $\Phi_n^\sigma$ is well defined linear map, provided
that we specify what ``finite sum'' means. So, finally, note that
the right hand side is the expression for the coefficient
$(a_{-1}b)_n$ in the $(-1)$-product \cite[Definition 3.7]{L2} of
mutually local twisted fields $Y^\sigma(a,z)$ and $Y^\sigma(b,z)$.
Hence there is certain $N(a,b)$ such that $a_jb=0$ for $j>N(a,b)$,
and we should take the last term to be
$$
\sum^{N}_{\substack{ l= 1}} (-1)^l\binom{\tfrac rT}{l}
\sum_{\substack{ j=-l \\i+j+1=m }}^{-1} \binom{l-1}{j+l}
[a_{j+\tfrac rT},b_{i+\tfrac sT}]
$$
for any $N\geq N(a,b)$. Our discussion shows that we have a linear map
\begin{equation*}
\Phi_n^\sigma \colon K_n^\sigma\rightarrow U^\sigma,\qquad
\Phi_n^\sigma \left(\sum_{i+j+1=n}a_i\otimes
b_j\right)=(a_{-1}b)_n,
\end{equation*}
where $(a_{-1}b)_n$ is the coefficient of the twisted vertex
operator $Y^\sigma(a_{-1}b,z)$.

Hence, by construction, we have the following consequence of
Proposition~\ref{T:Proposition 2.6,4}:

\begin{proposition}
\label{T:Proposition 2.6,5} For each $ q\in\ker \left(\Phi | \bar
R \mathbf 1\otimes N(k\Lambda_0)\right)$ and  $n\in\tfrac1T\mathbb
Z$ we have

\begin{equation}\label{2.14}
\Phi_n^\sigma(\chi_n^\sigma(q))=0.
\end{equation}
\end{proposition}

It is clear that $\Phi_n^\sigma(\chi_n^\sigma(q))=0$ is a relation
among coefficients of annihilating fields, but it will be more
convenient to say that $q_n=\chi_n^\sigma(q)$ itself is a relation
among coefficients of annihilating fields, or a relation among
relations. It will be also more convenient to consider
``coefficients $q_n$ of  vertex operators $Y^\sigma(q,z)$'' in a
slightly different way.

\subsection*{A connection between two formulations of relations among relations}

Let $a\in N(k\Lambda_0)$ be a homogeneous element of weight $\text{wt}\,  a$, i.e.,
$L_0a=(\text{wt}\, a) a$. Then the coefficient $a_i\in U^\sigma$ of the twisted vertex operator
$Y^\sigma(a,z)$ is of degree $(i+1-\text{wt}\, a)$, i.e., $[d,a_i]=-[L_0,a_i]=(i+1- \text{wt}\, a)a$.
Set $a_i=a(i+1-\text{wt}\, a)$. Then we have
$$
Y^\sigma(a,z)=\sum_{i\in\tfrac 1T\mathbb Z}a_iz^{-i-1}=
\sum_{n\in\tfrac 1T\mathbb Z}a(n)z^{-n-\text{wt}\, a}\,,
$$
where $a(n)\in U^\sigma(n)$. Likewise, for a homogeneous element
$q=a\otimes b$  in $N(k\Lambda_0)\otimes N(k\Lambda_0)$ of weight
$\text{wt}\, q=\text{wt}\, a+\text{wt}\, b$ we shall write
$$
Y^\sigma_{M\otimes M}(q,z)= \sum_{n\in\tfrac 1T\mathbb
Z}q(n)z^{-n-\text{wt}\, q}\,, \qquad q(n)=\sum_{i+j=n}a(i)\otimes
b(j)\,,
$$
where $q(n)$ are operators on $M\otimes M$. Again we want to make
sense of this formula for $a(i), b(j)\in U^\sigma$, but in a
slightly different way than before.

Set
$$
\big(U^\sigma\bar\otimes U^\sigma\big)(n)=
\prod_{i+j=n} \big(U^\sigma(i)\otimes U^\sigma(j)\big),
\qquad U^\sigma\bar\otimes U^\sigma=
\coprod_{n\in\tfrac 1T\mathbb Z} \big(U^\sigma\bar\otimes U^\sigma\big)(n).
$$
The elements of $U^\sigma\bar\otimes U^\sigma$ are linear
combinations of homogeneous sequences in $U^\sigma\otimes
U^\sigma$, usually we shall denote them as
$\sum_{i+j=n}a(i)\otimes b(j)$. For a fixed $n\in\frac1T\mathbb Z$
we have a linear map
$$
\chi^\sigma(n) \colon N(k\Lambda_0)\otimes
N(k\Lambda_0)\rightarrow \big(U^\sigma\bar\otimes U^\sigma\big)(n)
$$
defined for homogeneous elements $a$ and $b$ by
$$
\chi^\sigma(n) \colon a\otimes b\mapsto \sum_{p+r=n}a(p)\otimes
b(r).
$$
Again we think of $\chi^\sigma(n)(q)$ as ``the coefficient $q(n)$
of the vertex operator $Y^\sigma(q,z)$'', but now in another
space. We shall usually write $q(n)=\chi^\sigma(n)(q)$ for an
element $q\in N(k\Lambda_0)\otimes N(k\Lambda_0)$ and
$Q(n)=\chi^\sigma(n)(Q)$ for a subspace $Q\subset
N(k\Lambda_0)\otimes N(k\Lambda_0)$.

Since we have the adjoint action of $\hat{\mathfrak g}{[\sigma]}$
on $U^\sigma$, we define ``the adjoint action'' of $\hat{\mathfrak
g}{[\sigma]}$ on $U^\sigma\bar\otimes U^\sigma$ by
$$
[x(m),\sum_{p+r=n}a(p)\otimes b(r)]=
\sum_{p+r=n}[x(m),a(p)]\otimes b(r)+\sum_{p+r=n}a(p)\otimes
[x(m),b(r)].
$$
Note that we have the action of $\hat{\mathfrak g}$ on
$N(k\Lambda_0)\otimes N(k\Lambda_0)$ given by
$$
x_i(a\otimes b)= (x_i a)\otimes b+a\otimes (x_i b)\,,\qquad
x\in\mathfrak g, \ i\in\mathbb Z.
$$
As expected, we have the following commutator formula for
$q(n)=\chi^\sigma(n)(q)$:

\begin{proposition}
\label{T:Proposition 2.6,6} For\ $x(m)\in\hat{\mathfrak
g}{[\sigma]}$ and homogeneous $q$ we have
$$
[x(m),q(n)]=\sum_{i\geq 0}\binom{m}{i}(x_iq)(m+n),\qquad
(Dq)(n)=-(n+\text{wt}\, q)q(n).
$$
\end{proposition}

 So if a subspace $Q\subset N(k\Lambda_0)\otimes
N(k\Lambda_0)$ is invariant for $\coprod_{i\in\mathbb Z_{\geq
0}}{\mathfrak g}(i)$, then
$$
\coprod_{n\in\frac 1T\mathbb Z} Q(n)
$$
is a loop $\hat{\mathfrak g}{[\sigma]}$-module, in general
reducible.

Let us define the map $\Psi^\sigma \colon U^\sigma\bar\otimes
U^\sigma\rightarrow U^\sigma$ by
$$
\Psi^\sigma  \colon \sum_{p+r=n}a(p)\otimes b(r)\mapsto
\sum_{\substack{p+r=n \\ p<0}}a(p)b(r) +\sum_{\substack{p+r=n \\
p\geq 0}}b(r)a(p)
$$
(compare with the map $\Psi$ in  \cite{MP3}). Then for homogeneous $a\in N(k\Lambda_0)^r$, $b\in N(k\Lambda_0)^s$,
$n=m+\tfrac rT+\tfrac sT$, $n'=n+1-\text{wt}\, a-\text{wt}\, b$, we have
\begin{equation*}
\begin{split}
&\Phi_n^\sigma \big(\chi_n^\sigma(a\otimes b)\big)
- \Psi^\sigma \big(\chi^\sigma(n')(a\otimes b)\big)\\
&=-\sum^{\text{wt}\, a -2}_{\substack{ j=0\\i+j+1=m}} [a_{j+\tfrac
rT},b_{i+\tfrac sT}] +\sum_{\substack{ l\in\mathbb Z_{> 0}\\
\text{finite sum}}} (-1)^l\binom{\tfrac rT}{l} \sum_{\substack{
j=-l \\i+j+1=m }}^{-1} \binom{l-1}{j+l} [a_{j+\tfrac
rT},b_{i+\tfrac sT}].
\end{split}
\end{equation*}
By using this formula we can rewrite the relations (\ref{2.14})
among coefficients of annihilating fields in the following way
(compare with Proposition~3 in  \cite{MP3}):

\begin{proposition}
\label{T:Proposition 2.6,8} Let $q\in\ker \Phi$ be written as
\,$q=\sum a\otimes b$ in terms of homogeneous elements, and write
$r(a)=r\in\{0,1,\dots,T-1\}$ for $a\in N(k\Lambda_0)^r$. Then for
$n\in\tfrac 1T\mathbb Z$ we have
\begin{equation*}
\begin{split}
&\Psi^\sigma \big(q(n)\big)
 =\sum_{\substack{a,b}}\sum_{\substack{
i+j=n\\1-\text{wt}\, a \leq i\leq -1}} [a(i),b(j)]\\
&-\sum_{\substack{a, b}}\sum_{\substack{ l\in\mathbb Z_{> 0}\\
\text{finite sum}}} (-1)^l\binom{\tfrac {r(a)}{T}}{l}
\sum_{\substack{ i+j=n\\-l+1-\text{wt}\, a \leq i\leq -\text{wt}\,
a}} \binom{l-1}{1-\text{wt}\, a -i+\tfrac {r(a)}{T}} [a(i),b(j)].
\end{split}
\end{equation*}
\end{proposition}

Now assume that $q=\sum a\otimes b$ is a homogeneous element in
$\bar R{\mathbf 1}\otimes N(k\Lambda_0)$. Note that for $a\in\bar
R{\mathbf 1}$ the coefficient $a(i)$ of the corresponding field
$Y^\sigma(a,z)$ can be written as a finite linear combination of
basis elements $r(\rho)$, $\rho\in\ell \!\text{{\it t\,}}(\bar
R^\sigma)$. Hence each element of the sequence
$q(n)=\chi^\sigma(n)(\sum a\otimes b)\in \big(U^\sigma\bar\otimes
U^\sigma\big)(n)$, say $c_i$, can be written uniquely as a finite
sum of the form
$$
c_i=\sum_{\substack{ \rho\in\ell \!\text{{\it t\,}}(\bar R^\sigma)}} r(\rho)\otimes b_\rho,
$$
where $b_\rho\in U^\sigma$. If $b_\rho\neq 0$, then it is clear
that $|\rho|+|\ell \!\text{{\it t\,}}(b_\rho)|=n$. Let us assume
that $q(n)\neq 0$, and for nonzero ``$i$-th'' component $c_i$ let
$\pi_i$ be the smallest possible $\rho\ell \!\text{{\it
t\,}}(b_\rho)$ that appears in the expression for $c_i$. Denote by
$S$ the set of all such $\pi_i$. Since $q$ is a finite sum of
elements of the form $a\otimes b$, it is clear that there is
$\ell$ such that $\ell(\pi_i)\leq\ell$. Then, by our assumptions
on the order $\preceq$, the set $S$ has the minimal element, and
we call it the leading term $\ell \!\text{{\it
t\,}}\big(q(n)\big)$ of $q(n)$. For a subspace $Q\subset \bar
R{\mathbf 1}\otimes V$ set
$$
\ell \!\text{{\it t\,}}(Q(n))=\{\ell \!\text{{\it t\,}}(q(n))\mid
q\in Q,\, q(n)\neq 0\}.
$$
Note that our definition of leading terms depends on $\bar
R^\sigma$, and that we have not defined leading terms for general
elements in $\big(U^\sigma\bar\otimes U^\sigma\big)(n)$.

For a colored partition $\pi$ of set
$$
N(\pi)=\max \{\#\mathcal E(\pi)\!-\!1,\,0\},\quad
\mathcal E(\pi)=\{\rho\in \ell \!\text{{\it t\,}} (\bar R) \mid \rho\subset\pi\}.
$$

Note that $N(k\Lambda_0)\otimes N(k\Lambda_0)$ has a natural
filtration $\big(N(k\Lambda_0)\otimes N(k\Lambda_0)\big)_\ell$,
$\ell\in\mathbb Z_{\geq 0}$, inherited {}from the filtration
$\mathcal U_\ell$, $\ell\in\mathbb Z_{\geq 0}$.

\begin{lemma}[\cite{MP3}]
\label{T:Lemma 2.7} Let $Q\subset \ker \left(\Phi | (\bar R
\mathbf 1\otimes N(k\Lambda_0))_\ell\right)$
 be a finite dimensional subspace and $n\in\frac 1T\mathbb Z$.
Assume that $\ell(\pi)=\ell$ for all $\pi\in\ell \!\text{{\it
t\,}} (Q(n))$. Then
\begin{equation}\label{2.16}
\sum_{\pi\in\ell \!\text{{\it t\,}} (Q(n))} N(\pi)\geq \dim Q(n).
\end{equation}
Moreover, if in (\ref{2.16}) the equality holds, then for any two embeddings
 $\rho_1 \subset \pi$, $\rho_2 \subset \pi$, where $\rho_1, \rho_2 \in\ell \!\text{{\it t\,}}(\bar{R}^\sigma)$
and $\pi\in\ell \!\text{{\it t\,}} \big(Q(n)\big)$, we have a relation
\begin{equation}\label{2.17}
u(\rho_1 \subset \pi) \in  u(\rho_2 \subset \pi) + \overline
{\mathbb C\text{-span}\{u(\rho \subset \pi')\mid
 \rho \in \ell \!\text{{\it t\,}} (\bar{R}^\sigma), \rho \subset \pi', \pi \prec \pi'\}}.
\end{equation}
\end{lemma}
\begin{proof}
Let $Q_{[\pi]}=\{q\in Q(n) \mid \ell \!\text{{\it t\,}} (q) \succeq \pi\}$ and
$Q_{(\pi)}=\{q\in Q(n) \mid \ell \!\text{{\it t\,}} (q) \succ \pi\}$. Let
$\dim Q_{[\pi]}/Q_{(\pi)}=m(\pi)=m\leq \dim Q(n)$,
$m\geq 1$ and let $\rho_1\subset\pi,\dots ,\rho_s\subset\pi$ (where $s=s(\pi)$)
be all possible
embeddings in $\pi$. Let $\pi^*$ be such that $Q_{(\pi)}=Q_{[\pi^*]}$.
Let us follow the proof of Lemma 4 in \cite{MP3}:
we can write a basis of $Q_{[\pi]}/Q_{(\pi)}$ in the form
$$\aligned
c_{11} r(\rho_1)\otimes u(\pi/\rho_1)+&\dots+c_{1s} r(\rho_s)\otimes u(\pi/\rho_s)
+v_1+Q_{(\pi)},\\
c_{21} r(\rho_1)\otimes u(\pi/\rho_1)+&\dots+c_{2s} r(\rho_s)\otimes u(\pi/\rho_s)
+v_2+Q_{(\pi)},\\
&\dots\\
c_{m1} r(\rho_1)\otimes u(\pi/\rho_1)+&\dots+c_{ms} r(\rho_s)\otimes u(\pi/\rho_s)
+v_m+Q_{(\pi)},
\endaligned $$
where vectors $v_i$ are of the form
$$
v_i=\sum_{\pi\prec \pi'\prec \pi^*}\sum_{\rho\in\mathcal E(\pi')}
r(\rho)\otimes d_{i,\rho,\pi'}+
\sum_{\pi^*\preceq \pi'}\sum_{\rho\in\mathcal E(\pi')}
r(\rho)\otimes e_{i,\rho,\pi'} \in \big(U^\sigma\bar\otimes U^\sigma\big)(n),
$$
with $d_{i,\rho,\pi'}, e_{i,\rho,\pi'}\in U^\sigma$, $\ell \!\text{{\it t\,}}(d_{i,\rho,\pi'})=\pi'/\rho$,
$\ell \!\text{{\it t\,}}(e_{i,\rho,\pi'})=\pi'/\rho$.

Assume that $\text{rank\,}(c_{ij})<m$. Then the rows are linearly dependent.
By taking a nontrivial linear combination of the basis elements
we get a vector in $Q(n)$ of the form
$$
v=\sum_{\pi\prec \pi'\prec \pi^*}\sum_{\rho\in\mathcal E(\pi')}
r(\rho)\otimes d_{\rho,\pi'}+
\sum_{\pi^*\preceq \pi'}\sum_{\rho\in\mathcal E(\pi')}
r(\rho)\otimes e_{\rho,\pi'} .
$$
The elements $d_{\rho,\pi'}\in U^\sigma$, $\rho\in\mathcal E(\pi')$, $\pi\prec \pi'\prec \pi^*$,
must be zero since otherwise $Q_{(\pi)}\neq Q_{[\pi^*]}$. But then
$v\in Q_{(\pi)}=Q_{[\pi^*]}$ and our nontrivial linear combination of
basis elements is zero in $Q_{[\pi]}/Q_{(\pi)}$, a contradiction.

Since elements in $\big(U^\sigma\bar\otimes U^\sigma\big)(n)$ are really
the sequences, we should check that the above argument makes sense;
it works because the ``summands'' $r(\rho_1)\otimes u(\pi/\rho_1)$, for example, are always
at the same place in a sequence, depending on a degree $|\rho_1|$
of $r(\rho_1)$. (Note that the above argument would not make sense for
sequences in $\big(U^\sigma\otimes U^\sigma)^{\tfrac1T\mathbb Z}$.)

Hence the rank of matrix $(c_{ij})$ is $m$ and we have $s\geq m$.
Assume that $s=m$. Then the matrix $(c_{ij})$ is regular and the
vectors of the form $r(\rho_1)\otimes
u(\pi/\rho_1)+v_1+Q_{(\pi)}$,\dots , $r(\rho_m)\otimes
u(\pi/\rho_m)+v_m+Q_{(\pi)}$ (vectors $v_i$ as above) are a basis
of $Q_{[\pi]}/Q_{(\pi)}$. In particular, we have a vector $u\in
Q(n)$ of the form $u= r(\rho_1)\otimes u(\pi/\rho_1)+v_1$. {}From
the definition of $\Psi^\sigma$ we have
$$
\Psi^\sigma  (u)\in u(\pi)+U_{(\pi)}^\mathcal P.
$$
On the other hand $u=q(n)$ for some $q\in Q\subset \ker\Phi$, and
by assumption $\ell \!\text{{\it t\,}} (u)=\ell \!\text{{\it t\,}}
(q(n))=\ell$, so Proposition~\ref{T:Proposition 2.6,8} and
$\ell(\ell \!\text{{\it t\,}}(ab))>\ell(\ell \!\text{{\it
t\,}}([a,b]))$ imply
$$
\Psi^\sigma  (u)\in U_{(\pi)}^\mathcal P,
$$
what is in contradiction with Lemma~\ref{T:Lemma 2.2}.

Hence $s>m$, i.e. $N(\pi)= s(\pi)-1\geq m(\pi)$. Since
$\dim Q(n)=\sum m(\pi)$, the inequality (\ref{2.16}) follows.

Let us assume that $N(\pi)=m(\pi)$. Then there are altogether $s=m+1$ possible
embeddings $\rho_1\subset\pi,\dots ,\rho_s\subset\pi$ and the rank of $(c_{ij})$
is $m=s-1$. Let us assume that the first $s-1$ columns of $(c_{ij})$ are linearly
independent. Then for each $1\leq i\leq s-1$ there is a vector in $Q(n)$ of the form
$$
r(\rho_i)\otimes u(\pi/\rho_i)+d_i r(\rho_s)\otimes u(\pi/\rho_s)
+\sum_{\pi\prec \pi'}\sum_{\rho\in\mathcal E(\pi')}
r(\rho)\otimes d_{i,\rho,\pi'}
$$
for some $d_i\in \mathbb C$ and $d_{i,\rho,\pi'}\in U^\sigma$, $\ell \!\text{{\it t\,}}(d_{i,\rho,\pi'})=\pi'/\rho$,
As it was seen before, the assumption
$Q\subset ker \Phi$ implies $d_i\neq 0$. But then for each $i,j\in\{1,\dots ,s\}$
there are $d_{ij}\in \mathbb C$, $d_{i,j,\rho,\pi'}\in U^\sigma$ such that
$$
r(\rho_i)\otimes u(\pi/\rho_i)+d_{ij} r(\rho_j)\otimes u(\pi/\rho_j)
+\sum_{\pi\prec \pi'}\sum_{\rho\in\mathcal E(\pi')}
r(\rho)\otimes d_{i,j,\rho,\pi'}\in Q(n).
$$
By applying $\Psi^\sigma $ to these vectors, and taking into
account Proposition~\ref{T:Proposition 2.6,8}, we obtain relations
(\ref{2.17}).
\end{proof}

Even for low dimensional $Q$ it is not easy to determine the set
$\ell \!\text{{\it t\,}} (Q(n))$ of leading terms of $Q(n)$. For
this reason it is preferable to use the following simple
consequence of Lemma~\ref{T:Lemma 2.7}:

\begin{theorem}
\label{T:Theorem 2.8} Let $Q\subset \ker \left(\Phi | (\bar R
\mathbf 1\otimes N(k\Lambda_0))_\ell\right)$
 be a finite dimensional subspace and $n\in\frac 1T\mathbb Z$.
Assume that $\ell(\pi)=\ell$ for all $\pi\in\ell \!\text{{\it
t\,}} (Q(n))$. If
\begin{equation}\label{2.18}
\sum_{\pi\in\mathcal P^\ell(n)} N(\pi)=\dim Q(n),
\end{equation}
then for any two embeddings
 $\rho_1 \subset \pi$ and $\rho_2 \subset \pi$ in $\pi\in\mathcal P^\ell(n)$,
where $\rho_1, \rho_2 \in\ell \!\text{{\it t\,}}(\bar{R}^\sigma)$, we have a relation
\begin{equation}\label{2.19}
u(\rho_1 \subset \pi) \in u(\rho_2 \subset \pi) + \overline
{\mathbb C\text{-span}\{u(\rho \subset \pi')\mid
 \rho \in \ell \!\text{{\it t\,}} (\bar{R}^\sigma), \rho \subset \pi', \pi \prec \pi'\}}.
\end{equation}
\end{theorem}
\begin{proof}
Since our assumptions imply $\ell \!\text{{\it t\,}}
\big(Q(n)\big)\subset\mathcal P^\ell(n)$, we have
$$
\sum_{\pi\in\mathcal P^\ell(n)} N(\pi)\geq \sum_{\pi\in\ell
\!\text{{\it t\,}} (Q(n))} N(\pi).
$$
Hence our assumption (\ref{2.18}) implies the equality in
(\ref{2.16}), and the theorem follows {}from Lemma~\ref{T:Lemma
2.7}.
\end{proof}

\begin{remark}
By using the same arguments we can see that Theorem~\ref{T:Theorem
2.8} holds as well for $Q\subset \ker \left(\Phi |
(N(k\Lambda_0)\otimes \bar R \mathbf 1)_\ell\right)$.
\end{remark}

\begin{remarks}
(i) Note that the assumption (\ref{2.18}) relates purely
combinatorial quantity on the left hand side, which depends only
on $(\mathcal P,\preceq)$ and $\ell \!\text{{\it t\,}}
(\bar{R}^\sigma)$, with purely algebraic quantity on the right
hand side, which depends only on the kernel of the map $\Phi |
\bar R \mathbf 1\otimes N(k\Lambda_0)$.

It seems that ``the right hand side'' $\ker \left(\Phi | \bar R
\mathbf 1\otimes N(k\Lambda_0)\right)$ is easier to understand; it
can be described in terms of the kernel of the $\tilde{\mathfrak
g}$-module map
$$
 U(\tilde{\mathfrak g})\otimes_{U(\tilde{\mathfrak
g}_{\geq 0})}\bar R \mathbf 1 \rightarrow N(k\Lambda_0),\qquad
u\otimes v \mapsto uv
$$
between (induced) $\tilde{\mathfrak g}$-modules.

(ii) If $\ell(\rho)\leq K$ for all $\rho\in\ell \!\text{{\it t\,}}
(\bar{R}^\sigma)$, then, for constructing a combinatorial basis,
it is enough to construct combinatorial relations among relations
(\ref{2.19}) for all $\pi$ such that $\ell(\pi)\leq 2K-1$ (cf.
\cite{MP2, MP3}). In particular, for standard modules of level
$k$, and $R$ generated by the singular vector
$x_\theta(-1)^{k+1}\mathbf 1$, it is enough to construct
combinatorial relations among relations (\ref{2.19}) for all $\pi$
such that $\ell(\pi)\leq 2k+1$.

(iii) The condition that $\ell(\pi)=\ell$ for all $\pi\in\ell
\!\text{{\it t\,}} (Q(n))$ does not require a precise knowledge of
$\ell \!\text{{\it t\,}} (Q(n))$. For standard level $k$ modules,
and $R$ as above, this condition is obviously fulfilled for
$\ell=k+2$. In particular, for all level 1 standard modules this
condition is satisfied for $\ell=3$, and as observed above, it is
sufficient to consider this case alone. Hopefully, a precise
knowledge of highest weight vectors in $Q$, and the use of
Proposition~\ref{T:Proposition 2.6,6} with related loop module
structure, should make possible to check this condition for
general $\ell>k+2$.

(iv) A special case of Theorem~\ref{T:Theorem 2.8} for
$\mathfrak{g}=\mathfrak{sl}(2,\mathbb C)$ is implicit in the proof
of \cite[Lemma 9.2]{MP2}. In a general level $k\geq 2$ case the
relations (\ref{2.19}) for $\ell\leq 2k+1$ are obtained by
``explicitly solving'' certain systems of equations using the
relations for $\ell=k+2$. Having this example in mind, one may
think of Theorem~\ref{T:Theorem 2.8} as a kind of ``rank
theorem'', the ``existence of solutions'' (\ref{2.19}) follows
{}from calculating two sides of (\ref{2.18}). Of course, it may be
that in general finding the left hand side of (\ref{2.18}) is just
as hard as it is to ``explicitly solve'' the systems of equations
in question.

(v) Assume that $\sigma=\text{id}$ and that the chosen basis $B$
consists of weight vectors for the action of Cartan subalgebra
$\mathfrak h$. Then we can define $\mathfrak h$-weights of colored
partitions and Theorem~\ref{T:Theorem 2.8} can be refined in a
sense that for any $\mathfrak h$-weight $\mu$ the equality
$$
\sum_{\pi\in\mathcal P^\ell(n)_\mu} N(\pi)=\dim Q(n)_\mu
$$
implies (\ref{2.19}) for any two embeddings in $\pi\in\mathcal
P^\ell(n)_\mu$.

(vi) Theorem~\ref{T:Theorem 2.8} might be useful only if
(\ref{2.18}) holds in some generality. For now there are only
examples for twisted $A_1^{(1)}$ and $A_2^{(1)}$-modules. In the
case $\sigma=\text{id}$ for the basic $A_2^{(1)}$-module the
relation (\ref{2.18}) holds for all but two weights $\mu$ (in the
sense of previous remark, cf. \cite{MP3}). In the case when
$\sigma$ is a Dynkin diagram automorphism, there is a certain
basis $B$ of $\mathfrak{sl}(3,\mathbb C)$ such that for the
$\sigma$-twisted level 1 module (\ref{2.18}) holds (see \cite{S}).
In both examples the left hand side of (\ref{2.18}) is obtained by
direct counting of embeddings. For the right hand side the
appropriate $Q$ is guessed correctly, $\dim Q$ is calculated by
using the Weyl dimension formula, and then $\dim Q(n)$ is deduced
{}from the loop module structure related to $\coprod_{n\in\frac
1T\mathbb Z} Q(n)$.
\end{remarks}

\end{document}